\documentclass[oneside,a4paper]{article}
\usepackage{pp}
\usepackage{pictexwd,dcpic}
\usepackage[russian,english]{babel}
\usepackage[utf8]{inputenc}

\begin{document}

\title{Sweeping out sectional curvature}
\author{D. Panov\thanks{is a Royal Society University Research Fellow}
\ and A. Petrunin\thanks{was partially supported by the
Simons Foundation grant 245094.}}
\date{}

\maketitle

\begin{abstract}
We observe that the maximal open set of constant curvature $\kappa$
in a Riemannian manifold of curvature $\ge\kappa$ or $\le \kappa$
has a convexity type property, which we call \emph{two-convexity}.
This statement is used to prove a number of rigidity statements in
comparison geometry.
\end{abstract}

\maketitle

\section{Introduction}

Denote by $\MM^m[\kappa]$ the model $m$-space with curvature $\kappa$;
i.e., $\MM^m[\kappa]$ is the simply connected $m$-dimensional Riemannian manifold with constant curvature $\kappa$.
We will also use shortcuts $\mathbb S^m=\MM^m[1]$ for the unit $m$-sphere,
and $\EE^m=\MM^m[0]$ for the Euclidean $m$-space.

In this paper we play with applications of the following lemma.
The proof is given in Section~\ref{sec:key-lenmma}.
This lemma was first discovered by Buyalo
in the case of nonpositive curvature;
see \cite[Lemma 5.8]{buyalo}.

\begin{thm}{Buyalo's lemma}\label{lem:key}
Let $M$ be a complete Riemannian manifold with sectional curvature either $\ge \kappa$ or $\le \kappa$.
Let $\Delta$ be a tetrahedron in $\MM^3[\kappa]$
and $\Lambda$ be a union of three out of four faces of $\Delta$.
Then any immersion $f\:\Lambda\looparrowright M$ which is isometric and geodesic on each face
can be extended to an isometric  geodesic immersion $F\:\Delta\looparrowright M$.
Moreover, $F$ is uniquely determined by $f$.
\end{thm}

Here is an immediate corollary:

\begin{thm}{Corollary}
Let $g$ be a complete Riemannian metric on $\RR^3$ with curvature $\ge 0$ (or $\le 0$)
such that all three coordinate planes of $\RR^3$
are flat geodesic hypersurfaces in $(\RR^3,g)$.
Then $(\RR^3,g)$ is isometric to Euclidean space.
\end{thm}

We would suggest that reader checks that the last statement does not follow from the standard theorems;
in particular the splitting theorems cannot help here directly.

\medskip

Let us now introduce some terminology to state further applications.
\begin{itemize}
\item A Riemannian manifold (possibly not complete) of constant curvature $\kappa$ will be called \emph{$\kappa$-flat}.
\item A $\kappa$-flat Riemannian manifold (possibly not complete)
which satisfies the conclusion of Buyalo's lemma will be called \emph{two-convex}.
This definition is discussed in more detail in Section~\ref{sec:k-convex}.
\item Given a Riemannian manifold $M$, its maximal open subset of constant curvature $\kappa$ will be called the \emph{$\kappa$-flat domain of $M$}
and it will be denoted as $\Flat^\kappa M$.
\end{itemize}

From Buyalo's lemma one easily gets the following;
a formal proof is given in Section \ref{sec:key-lenmma}.

\begin{thm}{Observation}\label{obs:flat-domain}
Let $M$ be a complete Riemannian manifold either with curvature $\ge \kappa$ or $\le \kappa$.
Then $\Flat^\kappa M$ is two-convex.
\end{thm}

Here is an application.

\begin{thm}{Theorem}\label{thm:rigid-sphere}
Let $m\ge 3$ and $M$ be a complete connected $m$-dimensional manifold
with curvature $\ge 1$ or $\le 1$ which admits a totally geodesic immersion of the closed unit hemisphere $\iota\:\mathbb S^2_+\looparrowright M$
and an open neighborhood of $\iota(\mathbb S^2_+)$ in $M$ has constant curvature~$1$.
Then $M$ has constant curvature $1$.
\end{thm}

\parbf{Remarks.}
\begin{itemize}
\item Note that diameter-sphere rigidity does not help here directly;
in principle, the diameter of $M$ might be $<\pi$.

\item Note that  $\CP^2$ equipped with the canonical metric
is an example of a space with curvature $\ge 1$ and $\le 4$,
which admits totally geodesic immersions of 2-spheres of constant curvature $1$ and $4$.
I.e., the condition in Theorem~\ref{thm:rigid-sphere}
that the curvature is constant in a neighborhood of $\iota(\mathbb S^2_+)$ is necessary.

\item In the case of curvature $\ge 1$,
Theorem~\ref{thm:rigid-sphere} also holds in dimension $2$;
this is proved by Zalgaller in \cite{zalgaller};
see Theorem \ref{thm:rigid-sphere-2D} and the discussion around it.
\end{itemize}

\medskip

To prove the theorem,
one needs to show that if a neighborhood $\Omega$ of $\mathbb S^2_+$
in $\mathbb S^m$
admits an immersion in a two-convex manifold $\Phi$ then $\Phi$ has to be complete.
Then  Observation~\ref{obs:flat-domain} implies that $\Flat^1 M=M$;
i.e., $M$ is a spherical space form.
In other words, any neighborhood $\Omega$ of $\iota(\mathbb S^2_+)$ in $\mathbb S^m$ is \emph{exhaustive} in the sense of the following definition.

\begin{thm}{Definition}\label{def:exhaustive}
Let $\Omega$ be a $\kappa$-flat manifold.
Assume that any connected two-convex manifold $\Phi$ that appears as the target of an open isometric immersion $\Omega\looparrowright\Phi$ is complete
.
Then we say that $\Omega$ is \emph{exhaustive}.
\end{thm}

Using this definition,
we can formulate the following generalization of Theorem~\ref{thm:rigid-sphere}:

\begin{thm}{Theorem}\label{thm:exhaustive}
Let $M$ be a complete connected Riemannian manifold with curvature $\ge \kappa$ or $\le \kappa$.
Assume there is an open isometric immersion $\Omega\looparrowright M$
from an exhaustive $\kappa$-flat manifold $\Omega$.
Then $M$ has constant curvature $\kappa$.
\end{thm}

In order to apply this theorem one only has to find a source of exhaustive manifolds.
In Section~\ref{sec:k-convex}, we introduce the notion of the {\it two-hull} of a $\kappa$-flat simply connected  manifold $\Omega$;
in some sense this is the minimal simply connected two-convex manifold which contains  an \emph{immersed copy} of $\Omega$.
It is easy to see that  if the two-hull of a manifold $\Omega$ is isometric to $\MM^m[\kappa]$ then $\Omega$ is exhaustive.
This permits one to present a number of examples of exhaustive manifolds.
This is done in Section~\ref{sec:exhaustive},
here is a list of examples:

\medskip

\begin{description}
\item[(Proposition~\ref{prop:comp-conv}.)] For $m\ge 3$, any non-empty open subset of $\MM^m[\kappa]$ with convex complement.
\medskip

\item[(Proposition~\ref{prop:comp-conv+}.)]
More generally:
any open simply connected subset $\Omega\z\subset\MM^m[\kappa]$ which satisfies the following property.
For any $p\in \MM^m[\kappa]$
there is a 3-dimensional subspace $W_p$ of $\MM^m[\kappa]$ containing $p$
($W_p$ is an isometric copy of $\MM^3[\kappa]$)
such that $W_p\cap\Omega\not=\emptyset$
and each connected component of $W_p\backslash\Omega$ is a convex set.

In particular,
$$\Omega
=
\set{(x_1,x_2,\dots,x_m)\in \EE^m}{1+x_1^2+x_2^2>x_3^2+x_4^2+\dots+x_m^2}
$$
is exhaustive.

\medskip

\item[(Proposition~\ref{prop:S-m}.)] Any  open subset of $\mathbb S^m$ which
contains the standard 2-dimensional hemisphere.
This type of manifolds is used in Theorem~\ref{thm:rigid-sphere}.
\medskip
\end{description}
(This list can be continued.)

\parbf{Postscript.}
After publication of this note we learned that  Lars Andersson and Ralph Howard \cite{andersson-howard} obtained very similar results with a different method.

\subsection*{Related results.}

One outcome of Theorem~\ref{thm:exhaustive} is a sufficient condition on the closed set $K$%
\footnote{$K$ is the complement of $\Omega$} of the model space $\MM^m[\kappa]$,
such that one \emph{cannot} cut $K$ and glue instead a patch with sectional curvature either $\le \kappa$ or $\ge \kappa$ at all points. 
This condition is nontrivial only for $m\ge 3$.

The similar conditions for scalar and Ricci curvature were studied.
The case of deformation with nondecreasing  curvature turned out to be very different from the one with
nonincreasing curvature. 

After rescaling one can only consider three cases $\kappa=-1,0$ or $1$.

\parit{Nondecreasing curvature.}
If $\kappa=0$,
the case of nondecreasing scalar curvature leads to so-called
\emph{positive mass conjecture} which was proved by Schoen--Yau and 
Witten;
see
\cite{schoen-yau} and 
\cite{witten}.
This implies in particular that
\textit{the metric of Euclidean space cannot be perturbed in a bounded region so that the scalar curvature does not decrease}.

An analogous statement holds for $\kappa=-1$;
i.e., \textit{the metric of Lobachevsky space cannot be perturbed in a bounded region so that the scalar curvature does not decrease}.
The latter was proved by Min-Oo in \cite{oo-lob}.

The case $\kappa=1$ was considered in \cite{oo},
where Min-Oo  makes an attempt to show that the standard metric on the $m$-sphere
cannot be perturbed inside of hemisphere so that the scalar curvature does
not decrease.
Later, in \cite{brendle-marques-neves}, Brendle, Marques and Neves found a counterexample.
It is true that one cannot perturb the metric in a sufficiently small domain of the sphere,
but the optimal bounds for its size seem to be unknown.

The analogous statement for Ricci curvature was proved by Hang and Wang in \cite{hang-wang}.
They show that one cannot perturb the metric of the standard sphere inside its hemisphere with nondecreasing Ricci curvature.

The two-dimensional case of the above statements for $\kappa=0$ and $-1$ follows from Gauss--Bonnet formula
and the case $\kappa=1$ was done by Zalgaller (see the Appendix).

\parit{Nonincreasing curvature.}
In \cite{lohkamp-Ricc}, Lohkamp proves that
for all $m\ge3$,
one can perturb the metric of $\MM^m[\kappa]$ in any open region
in such a way that its Ricci curvature does not increase.
Moreover, this can be done without changing the topology 
and with an arbitrary small change to the geometry of the space.

In the two-dimensional case, attaching a handle can be done in an arbitrary small region with decreasing its curvature.
On the other hand, if we fix the topology, for $\kappa=0$ and $-1$,
Gauss--Bonnet formula prevents any change of metric in bounded regions
with nonincreasing curvature.
For $\kappa=1$, even if topology is fixed,
the metric can be changed (by inserting a bubble) in an arbitrary small open subset
so that the curvature in the region decreases. However, it seems that for proper
subsets of a hemisphere, there is no continuous deformation of this type.

\section{Two-convexity and two-hull}\label{sec:k-convex}

\begin{thm}{Definition}\label{def:2-conv-2-alt}
Let $\Omega$ be an $m$-dimensional $\kappa$-flat manifold.
We say that $\Omega$ is \emph{two-convex} if the following condition holds:
given a tetrahedron\footnote{I.e., 3-simplex.} $\Delta$ in $\MM^3[\kappa]$
with a choice of a subset $\Lambda\subset \Delta$ formed by 3 out of 4 faces,
any immersion $f\:\Lambda\looparrowright M$ which is isometric and geodesic on each face of $\Lambda$ can be extended to an isometric  geodesic immersion $F\:\Delta\looparrowright M$.
\end{thm}

\begin{thm}{Definition}\label{def:2-hull-alt}
Let $\Omega$ be a simply connected $m$-dimensional $\kappa$-flat manifold.
A simply connected $\kappa$-flat two-convex manifold $\Phi$
is called the \emph{two-hull} of $\Omega$ (briefly $\Phi=\Omega^{(2)}$)
if there is an open immersion $\phi\:\Omega\looparrowright\Phi$
such that for any open isometric immersion $\psi\:\Omega\z\looparrowright\Psi$ into a simply connected $\kappa$-flat two-convex manifold $\Psi$
there is a isometric immersion $\theta_\psi\:\Phi\z\looparrowright\Psi$ that makes the following diagram commutative:

$$\raisebox{-0.8cm}{$\begindc{\commdiag}[100]
\obj(3,5)[aa]{$\Omega$}
\obj(0,0)[bb]{$\Phi$}
\obj(6,0)[cc]{$\Psi$}
\mor{aa}{bb}{$\phi$}[-1,0]
\mor{aa}{cc}{$\psi$}
\mor{bb}{cc}{$\theta_\psi$}
\enddc$}
\eqno{({*})}
$$

In this case the immersion $\phi\:\Omega\looparrowright\Phi$ will be called the \emph{two-hull immersion}.
\end{thm}

Our next goal is to prove the existence of the two-hull.

\begin{thm}{Proposition}\label{prop:2-hull-alt}
For any simply connected $\kappa$-flat manifold $\Omega$,
its two-hull $\Phi$ is uniquely defined up to an isometry.

Moreover, if $\phi\:\Omega\looparrowright\Phi$ and $\phi'\:\Omega\looparrowright\Phi'$
are two-hull immersions
then there is an isometry $\theta:\Phi\to \Phi'$ such that
$\phi'=\theta\circ\phi$.
\end{thm}

To prove the above proposition, we mimic the proof of the existence of
ordinary convex hull obtained by intersecting all convex sets containing the given set.

\parit{Proof.}
Assume $\Psi$ is a simply connected $\kappa$-flat two-convex manifold and $\psi\:\Omega\z\looparrowright\Psi$ is an isometric immersion.

Since $\Omega$ is simply connected, it admits an isometric immersion $\iota:\Omega\looparrowright\MM^m[\kappa]$.
Moreover,
there is an immersion
$\tau_\psi\:\Psi\to \MM^m[\kappa]$ which makes the following diagram commutative.
$$\raisebox{-0.9cm}{$\begindc{\commdiag}[100]
\obj(3,5)[aa]{$\Omega$}
\obj(0,0)[bb]{$\Psi$}
\obj(6,0)[cc]{$\MM^m[\kappa]$}
\mor{aa}{bb}{$\psi$}[-1,0]
\mor{aa}{cc}{$\iota$}
\mor{bb}{cc}{$\tau_\psi$}
\enddc$}
$$

Fix a point $x\in\Omega$.
For any $\psi$ as above
set $x_\psi=\psi(x)\in \Psi$;
in particular, $x_\iota=\iota(x)\z\in \MM^m[\kappa]$.

Consider the set $\Gamma$ of all paths in $\MM[\kappa]$ which start at $x_\iota$.
Let us equip $\Gamma$ with $C^0$ topology.
Given $\psi$ as above,
we say that $\gamma\in \Gamma$ lifts to $\psi$
if there is a path $\gamma_\psi\:[0,1]\to\Psi$ which starts at $x_\psi$ and such that $\gamma=\tau_\psi\circ\gamma_\psi$.
The set of all paths which lift to $\psi$ will be denoted as $\Gamma_\psi$.

Assume $\Psi'$ is some other simply connected $\kappa$-flat two-convex manifold and $\psi'\:\Omega\z\looparrowright\Psi'$ is an isometric immersion.
Note that if $\Gamma_\psi\subset\Gamma_{\psi'}$ 
then there is a necessarily unique isometric immersion $\lambda\:\Psi\to \Psi'$
which makes the following diagram commutative:
$$\raisebox{-0.9cm}{$\begindc{\commdiag}[100]
\obj(3,5)[aa]{$\Omega$}
\obj(0,0)[bb]{$\Psi$}
\obj(6,0)[cc]{$\Psi'$}
\mor{aa}{bb}{$\psi$}[-1,0]
\mor{aa}{cc}{$\psi'$}
\mor{bb}{cc}{$\lambda$}
\enddc$}
$$
In particular if $\Gamma_\psi=\Gamma_{\psi'}$ 
then  $\lambda$ is an isometry.

Denote by $\Xi\subset\Gamma$ the interior of the space $\bigcap_\psi\Gamma_\psi$,
where the intersection is taken for all $\psi$ as above.

Denote by $\Phi$ the set of the homotopy classes of paths in $\Xi$ rel{.} the ends;
i.e., we consider only the homotopies $\gamma_t$ such that $\gamma_t\in \Xi$ for any $t$.
Denote the projection by $\pi\:\Xi\to\Phi$.
The set $\Phi$ comes with a topology and a $\kappa$-flat metric
which makes the maps $\theta_\psi\:\pi(\gamma)\mapsto \gamma_\psi(1)$ isometric immersions  $\theta_\psi\:\Phi\looparrowright\Psi$ 
for all $\psi$ as above. 

Given a point $y\in \Omega$ consider a path $\beta$ from $x$ to $y$ in $\Omega$.
Note that for any $\psi$ as above,
the path $\psi\circ\beta$ is a lift of $\iota\circ\beta$ to $\psi$;
the same holds for any path $\beta'\in\Gamma$ 
which is sufficiently close to $\beta$.
It follows that $\iota\circ\beta\in\Xi$.

Consider the map $\phi\:\Omega\to \Phi$, defined as $\phi(y)= \pi(\iota\circ\beta)$.
Note that the value $\phi(y)$ does not depend on the choice of $\beta$ since $\Omega$ is simply connected.
Moreover $\phi$ is an isometric immersion and
it makes the diagram $({*})$ commutative.

Summarizing all the above, $\phi\:\Omega\looparrowright\Phi$ is the two-hull immersion.

Finally, assume $\phi\:\Omega\looparrowright\Phi$ and $\phi'\:\Omega\looparrowright\Phi'$
are two two-hull immersions.
Then the immersions  $\theta\:\Phi\looparrowright\Phi'$ and $\theta'\:\Phi'\looparrowright\Phi$
provided by Definition~\ref{def:2-hull-alt} 
have to be inverses of each other.
The latter implies the last statement of the proposition.
\qeds

\section{Buyalo's lemma and the observation}\label{sec:key-lenmma}

In this section we prove Buyalo's lemma \ref{lem:key}
and Observation~\ref{obs:flat-domain}.
The proof of the following proposition is left to the reader.

\begin{thm}{Proposition}\label{lem:triv}
Let $X$ and $Y$ be (possibly noncomplete) Riemannian manifolds
and $\Gamma$ be an open set of unit speed geodesics in $X$,
covering all points of $X$.
Then $f\:X\to Y$ is an isometric geodesic immersion
if and only if for any $\gamma\in\Gamma$,
the curve $f\circ\gamma$ is a unit speed geodesic in $Y$.
\end{thm}

\parit{Proof of Buyalo's lemma.}
Set $m=\dim M$.
Note that the statement of Buyalo's lemma trivially holds if $m\le 2$.
Further we assume $m\ge 3$.

By choosing an isometric geodesic embedding $\Delta\hookrightarrow\MM^m[\kappa]$,
we can consider $\Delta$ as a subset of $\MM^m[\kappa]$.
Let us denote by $\tilde p$ the common vertex of the faces in $\Lambda$ and
let $\tilde x,\tilde y,\tilde z$ be the remaining vertexes of $\Delta$.
Denote by $p,x,y,z$ the corresponding points in $M$;
i.e.
\begin{align*}
p&=f(\tilde p),
&
x&=f(\tilde x),
&
y&=f(\tilde y),
&
z&=f(\tilde z).
\end{align*}

Fix $R>\diam \Delta$.
Assume first that the injectivity radius at any point in $B_{2\cdot R}(p)\subset M$ is at least $2\cdot R$.
In this case $f$ is  distance preserving on each face.

\begin{thm}{Claim}\label{clm:dist-preserv}
$f\:\Lambda\to M$ is a distance preserving map.
\end{thm}

\begin{wrapfigure}{r}{40mm}\begin{lpic}[t(-6mm),b(0mm),r(0mm),l(0mm)]{pics/pxyz(0.8)}\lbl[r]{23,28;$p$}\lbl[t]{13,0;$x$}\lbl[b]{2,45;$y$}\lbl[b]{49,45;$z$}\lbl[r]{17,14;{\small $x'$}}\lbl[b]{14,38;{\small $y'$}}\lbl[b]{31,34.5;{\small $z'$}}\lbl[r]{13.5,28;$v$}\lbl[b]{22,36;$w$}\end{lpic}\end{wrapfigure}

\parit{Proof of the claim.}
On the geodesic $[px]$ consider two unit normal fields
that go in the directions of the images of the faces adjacent to $[px]$.
Note that both fields are parallel.
Thus, the angle between the images of the faces in $\Lambda$ is constant along the common side.
Taking a point on the geodesic $[px]$ close to $p$,
one can see that the angles between faces of $f(\Lambda)$ in $M$
coincide with the corresponding angles in $\Lambda\subset \MM^m[\kappa]$.

Consider points
\begin{align*}
\tilde x'&\in [\tilde p\tilde x],
&
\tilde y'&\in [\tilde p\tilde y],
&
\tilde z'&\in [\tilde p\tilde z],
\\
x'&=f(\tilde x'),
&
y'&=f(\tilde y'),
&
z'&=f(\tilde z').
\end{align*}
From above, we have that corresponding angles in the triangles $[x'y'z']$ and $[\tilde x'\tilde y'\tilde z']$
are equal;
i.e., the angles in the triangle $[x'y'z']$ coincide with its comparison angles.

Let $\tilde v$ and $\tilde w$ be arbitrary points on the sides of the triangle $[\tilde x'\tilde y'\tilde z']$
and $v=f(\tilde v)$ and $w=f(\tilde w)$.
In both cases (curvature $\ge \kappa$ or $\le \kappa$)
angle-sidelength  monotonicity (see for example \cite{AKP}) implies that
$$|v-w|_M=|\tilde v-\tilde w|_{\MM^m[\kappa]},$$
where $|{*}-{*}|_X$ denotes the distance function in the metric space $X$.

Note that for any $\tilde v,\tilde w\in \Lambda$ there is a triangle $[\tilde x'\tilde y'\tilde z']$ as above
which contains $\tilde v$ and $\tilde w$ on its sides.
Hence, the claim follows.
\qeds

Note that there is a map  $F\:B_R(\tilde p)\to B_R(p)$ satisfying
the following properties:
\begin{enumerate}
\item $F|_\Lambda=f$;
\item $F(\tilde p)=p$,
and the differential of $F$ at $\tilde p$ is an isometry $\T_{\tilde p}\to \T_p$;
\item $F$ sends all unit speed geodesics through $\tilde p$ to unit speed geodesics through $p$.
\end{enumerate}

\begin{thm}{Claim}
The restriction of any such $F$ to $\Delta$ satisfies Buyalo's lemma;
\end{thm}

This claim is proved separately in the following two cases:

\parit{Proof of the claim in the case of curvature $\ge \kappa$.}
By Rauch comparison (see for example \cite[Corollary 1.35]{cheeger-ebin})
the diffeomorphism
$F\:B_R(\tilde p)\to B_R(p)$
is non-expanding. 
Together with Claim~\ref{clm:dist-preserv},
this implies that
 the restriction of $F$ to $\Delta$ is distance preserving on any geodesic in $\Delta$ with ends in $\Lambda$. 
 
 Applying Proposition~\ref{lem:triv}
we get that the restriction of $F$ to $\Delta$
is isometric and geodesic in the interior of $\Delta$
and hence the same holds on whole $\Delta$.\qeds

\parit{Proof of the claim in the case of curvature $\le \kappa$.}
Set $\Upsilon$ to be the set of all minimizing geodesics with ends in $f(\Lambda)$
and let $\bar \Upsilon$ be the subset of $M$ covered by all geodesics in $\Upsilon$.

By the Rauch comparison, the diffeomorphism $F\:B_R(\tilde p)\to B_R(p)$ is non-contracting,
while its inverse $F^{-1}$ is a non-expanding diffeomorphism.
Since $f$ is distance preserving, it follows that $F^{-1}$ is isometric on each of the geodesic in $\Upsilon$;
moreover, any minimizing geodesic between points in $\Lambda$ can be presented as $F^{-1}\circ\gamma$ for some $\gamma\in \Upsilon$.
It follows that $F^{-1}(\bar \Upsilon)=\Delta$, or equivalently $F(\Delta)=\bar\Upsilon$.
In particular, $F$ is distance preserving on each minimizing geodesic with ends in $\Lambda$.

Applying Proposition~\ref{lem:triv} the same way as above,
we conclude that  the restriction of $F$ to $\Delta$
is distance preserving and geodesic.\qeds

\parit{The general case.}
To treat the general case, choose $\eps>0$ so that the injectivity radius at any point in $B_{2\cdot R}(p)$ is at least $2\cdot\eps$.
Note that one can cover the interior of $\Delta$
by an infinite sequence of tetrahedra $\Delta_1,\Delta_2,\dots$
with a choice of three faces $\Lambda_i$ in each $\Delta_i$
such that $\diam \Delta_i<\eps$ and
$$\Lambda_n
\subset
\Lambda\cup\left(\bigcup_{i<n}\Delta_i\right).$$
for each $n$.
Then it remains to apply the above argument sequentially to $\Delta_1,\Delta_2,\dots$ and pass to the closure.
\qeds

\parit{Proof of Observation~\ref{obs:flat-domain}.}
Set $m=\dim M$.
Choose any point $p\in\Flat^\kappa M$ and $\tilde p \in \MM^m[\kappa]$.
Choose a map $F\:\MM^m[\kappa]\to \Flat^\kappa M$ such that
\begin{enumerate}
\item $F(\tilde p)=p$,
and the differential of $F$ at $p$ is an isometry $\T_{\tilde p}\to \T_p$;
\item $F$ sends all unit speed geodesics through $\tilde p$ to unit speed geodesics through $p$.
\end{enumerate}

Let $\Omega_{p}\subset \MM^m[\kappa]$ be the maximal open star-shaped w.r.t. $\tilde p$ set
such that the map $F$ induces an open isometric immersion of $\Omega_p$.
Let $\Psi_{p}$ be the set of all tetrahedra with one vertex at $\tilde p$ and three adjacent faces in $\Omega_{p}$
and let $\bar\Psi_{p}$ be the union of all tetrahedra in $\Psi_p$.

Clearly $\bar\Psi_{p}$ is open and $\bar\Psi_{p}\supset \Omega_{p}$.
According to Buyalo's lemma, the map $F$ is isometric on each geodesic lying in a tetrahedron from $\Psi_{p}$.
Applying Proposition~\ref{lem:triv}, we get that  $F$
is an open isometric immersion $\bar\Psi_{p} \looparrowright M$.
Thus, $\bar\Psi_{p}=\Omega_{p}$ for any $p\in \Flat^\kappa M$,
hence the result.
\qeds

\section{Exhaustive manifolds}\label{sec:exhaustive}

Let $\Omega$ be a simply connected $\kappa$-flat manifold.
Recall that $\Omega^{(2)}$ denotes the two-hull of $\Omega$ (see Definition~\ref{def:2-hull-alt}).
From the definition of the two-hull, we have that if $\Omega^{(2)}$ is isometric to the model space $\MM^m[\kappa]$ then $\Omega$ is exhaustive (see Definition~\ref{def:exhaustive}).

In this section we use the above observation to construct examples of exhaustive manifolds.
The following two propositions follow directly from the discussion above.
(In other words, the proof is left to the reader.)

\begin{thm}{Proposition}\label{prop:comp-conv}
Assume $m\ge 3$ and suppose $\Omega\subset\MM^m[\kappa]$ is a nonempty open set with convex complement. Then $\Omega^{(2)}$ is isometric to $\MM^m[\kappa]$.
In particular, $\Omega$ is exhaustive.
\end{thm}

Here is a generalization of the above proposition:

\begin{thm}{Proposition}\label{prop:comp-conv+}
Suppose $m\ge 3$ and suppose $\Omega\subset\MM^m[\kappa]$ is a nonempty open set such that through any point $p\in \MM^m[\kappa]$ passes a $3$-dimensional subspace $W_p$ (i.e., an isometric copy of $\MM^3[\kappa]$) such that each connected component of $W_p\backslash\Omega$ is a convex set.

Then $\Omega^{(2)}$ is isometric to $\MM^m[\kappa]$.
In particular, $\Omega$ is exhaustive.
\end{thm}

The proof of the following proposition requires some work.
Set $\mathbb S^m\df\MM^m[1]$.

\begin{thm}{Proposition}\label{prop:S-m}
Assume $m\ge 3$ and suppose that an open set 
$\Omega\subset\mathbb S^m$
admits a geodesic isometric immersion $\mathbb S^2_+\hookrightarrow \Omega$.
Then $\Omega^{(2)}$ is isometric to $\mathbb S^m$.
\end{thm}

\parit{Proof.}
Fix two embeddings $\mathbb S^2_+\hookrightarrow\Omega$ and $\Omega\hookrightarrow\mathbb S^m$,
denote their composition by $\iota$.
Note that for any point $x\in \mathbb S^m\backslash \iota(\partial\mathbb S^2_+)$
there is the unique embedding $\iota_x\:\mathbb S^2_+\hookrightarrow\mathbb S^m$
such that $x\in \iota_x(\mathbb S^2_+)$
and $\iota_x(z)=\iota(z)$ for any $z\in\partial\mathbb S^2_+$.
It is easy to see that one can choose a tetrahedron $\Delta$ in $\mathbb S^m$, such that one face of
$\Delta$ belongs to $\iota_x(\mathbb S^2_+)$ and contains all points in the set  $\iota_x(\mathbb S^2_+)\backslash\Omega$, while the rest of the faces are arbitrary close to $\iota(\mathbb S^2_+)$, in particular these faces belong to $\Omega$.

Applying to $\Delta$ the definition of two-convexity,
we get an isometric geodesic immersion $F\:\Delta\looparrowright\Omega^{(2)}$.
It is easy to see that the map $x\mapsto F(x)$ is independent of the choice of $\Delta$;
moreover, the obtained map $\mathbb S^m\to \Omega^{(2)}$ is an open isometric immersion.
Since $\Omega^{(2)}$ is simply connected,
we have that $\Omega^{(2)}$ is isometric to $\mathbb S^m$.
\qeds

\section{Comments and open problems}

\parbf{$\bm{k}$-convexity.}
The definition of two-convexity (\ref{def:2-conv-2-alt}) can be generalized to  ``$k$-convexity'';
one has to change the tetrahedron $\Delta$ to a $(k+1)$-dimensional simplex
and $\Lambda$ to the set formed by $k+1$ faces out of $k+2$ in $\Delta$.
In this case, $1$-convexity is equivalent to the usual convexity of each connected component of $\Omega$.

In \cite[Section~$\tfrac12$]{gromov},
Gromov introduced the following closely related notion which we will call further as \emph{Lefschetz-$k$-convexity}%
\footnote{We state a slight variation of Gromov's definition;
in particular, we restrict our consideration to open sets
and change the meaning of $k$; in Gromov's notations Lefschetz-$k$-convexity in $\EE^m$ is called $(m-k)$-convexity.}.

\begin{thm}{Definition}\label{def:lefschetz}
An open set $\Omega$ in $\EE^m$ is \emph{Lefschetz-$k$-convex}
if for any
$k$-dimensional affine subspace $A$ the natural homology homomorphism
$$H_{k-1}(\Omega\cap A)\to H_{k-1}(\Omega)\eqno({*}{*})$$
is injective.
\end{thm}

This definition can be generalized to $\kappa$-flat manifolds,
one only has to replace $\Omega\cap A$ by $k$-dimensional manifolds $\Theta$ that admit a proper%
\footnote{An isometric immersion $\iota\:\Theta\looparrowright\Omega$ of Riemannian manifolds $\Theta$ and $\Omega$ is called \emph{proper}
if for any point $p\in\Omega$ there is $\eps>0$ such that each connected component of $\iota^{-1}(\bar B_\eps(p))\subset \Theta$ is compact.} isometric geodesic immersion $\Theta\looparrowright\Omega$.

It is easy to show that Lefschetz-$k$-convexity in $\EE^m$ implies our $k$-convexity.
We know that the converse holds in two trivial cases: $k=1$ and $m\le k+1$,
but in all other cases we do not know the answer to the following question.

\begin{thm}{Open problem}
Is it true that any $k$-convex open subset of $\EE^m$ is Lefschetz-$k$-convex?
\end{thm}

\parbf{Smooth approximation of two-convex sets.}
To get a feeling of the definition of $k$-convexity,
it is useful to observe the following.

\begin{thm}{Proposition}\label{prop:smooth}
If $\Omega$ is an open subset of $\EE^m$ with smooth boundary $\partial\Omega$,
then it is $k$-convex if and only if the hypersurface
$\partial\Omega$ has at most $k-1$ negative principal curvatures at any point.
\end{thm}

It is well known that any convex set in $\EE^m$ can be approximated
 by a convex set with smooth boundary.
It turns out that for $k$-convex sets (as well as for Lefschetz-$k$-convex sets)
this is no longer true.

To give an explicit example, let $\Omega\subset\RR^4$ be the complement to the union of the following two 3-dimensional halfspaces: $$(x_1\ge ax_2, \, x_3=0)\cup(x_1\ge -ax_2,\,x_4=0),\,\, a>0$$
The intersection of the halfspaces lies in the two-plane $x_3=x_4=0$ and 
forms a plane angle $Q$ in it.
Clearly $\Omega$ is connected, 2-convex and simply connected.
Let us give two alternative ways to prove that $\Omega$ is not smoothable.

\parit{Way 1.}
Note that for $k$-convex sets with smooth boundary the homeomorphism in $({*}{*})$ is injective for subspaces $A$ of arbitrary dimension;
the proof is an exercise in Morse theory, see \cite[Section~$\tfrac12$]{gromov}.
Thus, any $k$-convex set which does not satisfy this condition cannot be approximated.

Let $A$ be the 3-subspace in $\EE^4$ given by the equation 
$x_1=0$.
Note that $A\cap\Omega$ is formed by the complement of
two half-planes in $A$ intersecting at single point.
In particular $H_1(A\cap\Omega)=\ZZ$, 
which can be shown to contradict the above observation.

\parit{Way 2.}
Assume $\Omega$ can be approximated
 by a $k$-convex set with smooth boundary.
Equip $\Omega$ with the intrinsic metric induced in $\EE^4$
and denote by $\hat\Omega$ its completion. 
According to the main result in \cite{ABB}, 
$\hat\Omega$ is $\CAT(0)$. 
The latter is not the case, say if the angle measure of $Q$ is less than $\tfrac\pi2$.

\parbf{Two-hull in the non simply connected case.}
The following example shows a problem with extension of the
two-hull construction to the  non simply connected case.
Consider an isometric action $\ZZ_2\acts\mathbb S^3$ with two fixed points;
then take $\Omega$ to be the orbit space $\mathbb S^3/\ZZ_2$ with singular orbits removed.
Note that $\Omega$ admits no open isometric immersion into a two-convex $1$-flat manifold.
Hence the two-hull of $\Omega$ cannot be defined in the class of manifolds.

Note that if a $\kappa$-flat manifold $\Omega$ does not admit a two-hull then
it is automatically exhaustive. In this case 
$\Omega$ is not isometric to a $\kappa$-flat subset in any manifold of curvature $\ge \kappa$ or $\le \kappa$.

On the other hand, the two-hull is always defined in the class of so-called \emph{Riemannian megafolds} (in particular Riemannian orbifolds);
these creatures were introduced in \cite{petrunin-tuschmann} and under a different name in  \cite{lott};
they look a lot like Riemannian manifolds,
but fail to be topological spaces.
In the above example, the two-hull of $\Omega$ is the Riemannian orbifold $(\mathbb S^3:\ZZ_2)$. 

\parbf{More questions.}
Here is a possible generalization of Proposition~\ref{prop:S-m}:

\begin{thm}{Question}
 Is it true that the two-hull of any open simply connected set $\Omega\subset\mathbb S^m$ which contains a closed geodesic is isometric to $\mathbb S^m$?
\end{thm}

The following question of Burago and Kleiner has been open for a long time.
It is not directly relevant to all above,
but it was one of the initial motivations for our work.

\begin{thm}{Question}
Is it possible to construct a Riemannian metric $g$ on the product of a
torus and an open disc $T^2\times D^2$
such that the torus $T^2\times \{0\}\z\hookrightarrow T^2\times D^2$ is flat
and the curvature is strictly positive outside $T^2\times \{0\}$?
\end{thm}

An answer to this question might lead to a better understanding of manifolds with \emph{almost positive curvature} (see \cite{ziller}).

\medskip

We also mention two related questions from \href{http://mathoverflow.net/}{mathoverflow}:
\begin{itemize}
 \item \href{http://mathoverflow.net/questions/55788}{Question 55788}
about two-convexity and the Lefschetz property.
\item \href{http://mathoverflow.net/questions/50889}{Question  50889}
about possible generalization of Buyalo's lemma.
\end{itemize}

\def\claim#1{\par\medskip\noindent\refstepcounter{thm}\hbox{\bf \Alph{section}.\arabic{thm}. #1.}
\it\ 
}

\appendix
\setcounter{section}{1}
\setcounter{thm}{0}
\section*{Appendix: Zalgaller's rigidity.}\label{subsec:zalgaller}
Here we briefly repeat the proof of a theorem from \cite{zalgaller}.
We do this since the result which interests us (Theorem~\ref{thm:rigid-sphere-2D}) was not formulated as a separate statement;
it appeared as an intermediate statement in the proof.

\begin{thm}{Theorem}\label{thm:zalgaller}
Let $A=a_1a_2\dots a_n$ and $B=b_1b_2\dots b_n$ be two simple spherical polygons (not necessary convex)
with equal corresponding sides.
Assume $A$ lies in an open hemisphere and $\angle a_i\ge\angle b_i$ for each $i$.
Then $A$ is congruent to $B$.

\end{thm}

At first this result might look unrelated to the content of this article.
But the proof relies on the following 2-dimensional analog of Theorem~\ref{thm:rigid-sphere}.
Recall that a \emph{spherical polyhedron} is a simplicial complex equipped with a metric such that each simplex is isometric to a simplex in a standard sphere.

\begin{thm}{Theorem}\label{thm:rigid-sphere-2D}
Let $\Sigma$ be a spherical polyhedron which is homeomorphic to $\mathbb S^2$ and has curvature $\ge 1$ in the sense of Alexandrov. Assume that an open neighborhood of $\mathbb S^2_+$ in $\mathbb S^2$
admits a locally isometric immersion in $\Sigma$.
Then $\Sigma$ is isometric to the standard sphere.
\end{thm}

To deduce Theorem~\ref{thm:zalgaller} from Theorem~\ref{thm:rigid-sphere-2D},
Zalgaller cuts the polygon $A$ from the sphere and glues the polygon $B$ there instead.
As a result he gets the spherical polyhedron $\Sigma$ as in Theorem~\ref{thm:rigid-sphere-2D}.
(In fact, if we drop the condition that $A$ lies in a hemisphere,
we can obtain in this way any spherical polyhedral metric on $\mathbb S^2$ with curvature $\ge 1$.)

Theorem~\ref{thm:rigid-sphere-2D} is proved by induction
on the number $n$ of singular points in $\Sigma$.
The base case $n=1$ is trivial.
To do the induction step,
choose two singular points $p, q\in\Sigma$,
cut $\Sigma$ along a geodesic $[pq]$
and patch the hole so that the obtained new polyhedron $\Sigma'$ has curvature $\ge 1$.
The patch is obtained by gluing two copies of a
spherical triangle along two sides.
For the right choice of the triangle,
the points $p$ and $q$ become regular in $\Sigma'$
and exactly one new singular point appears in the patch.
This way, the case with $n$ singular points is reduced
to the case with $n-1$ singular points (if $n>1$).

The patch construction above was introduced by Alexandrov
in his famous proof of convex embeddability of polyhedrons;
the earliest reference we have found is
\cite[VI, \S7]{alexandrov1948}.

Applying polyhedral approximation, one can extend Theorem~\ref{thm:rigid-sphere-2D}
to any surface with curvature $\ge 1$ in the sense of Alexandrov;
in particular, this shows that Theorem~\ref{thm:rigid-sphere} holds in addition for $m=2$ and curvature $\ge 1$.

\end{document}